\input ./style/arxiv-general.cfg
\documentclass[aos,MSNbibl,nameyear,dvips]{arximspdf}
\makeatletter
   \@ifpackageloaded{graphicx}{}{\usepackage{graphicx}}
\makeatother

\doi{10.1214/15-AOS1359}
\volume{43}
\issue{6}
\pubyear{2015}
\firstpage{2738}
\lastpage{2765}
\docsubty{FLA}

\makeatletter
\newcommand{\rrvert}{\vert}
\newcommand{\llvert}{\vert}
\newproclaim{rem}{Remark}

\def\b{\beta}

\def\oo{\infty}

\def\u{\mathbf{u}}
\def\x{\mathbf{x}}
\def\b{\mathbf{b}}

\def\X{\mathbf{X}}
\def\A{\mathbf{A}}
\def\Y{\mathbf{Y}}

\def\vv{\mathbf{v}}
\def\1{\mathbf{1}}

\newtheorem{theorem}{Theorem}
\newtheorem{lemma}{Lemma}
\makeatother

\begin{document}
\begin{frontmatter}

\title{Bridging centrality and extremity:
Refining empirical data depth using extreme value statistics}
\runtitle{Data depth and extreme values}

\begin{aug}
\author[A]{\fnms{John H.~J.}~\snm{Einmahl}\corref{}\ead[label=e1]{j.h.j.einmahl@uvt.nl}},
\author[B]{\fnms{Jun}~\snm{Li}\ead[label=e2]{jun.li@ucr.edu}}
\and
\author[C]{\fnms{Regina Y.}~\snm{Liu}\thanksref{T1}\ead[label=e3]{rliu@stat.rutgers.edu}}
\runauthor{J. H.~J. Einmahl, J. Li and R.~Y. Liu}
\thankstext{T1}{Supported by NSF Grants DMS-07-07053 and DMS-10-07683.}
\affiliation{Tilburg University, University of California, Riverside\break
and Rutgers University}
\address[A]{J. H.~J. Einmahl\\
Department of Econometrics and CentER\\
Tilburg University\\
P.O. Box 90153\\
5000 LE Tilburg\\
The Netherlands\\
\printead{e1}}
\address[B]{J. Li\\
Department of Statistics\\
University of California, Riverside\\
Riverside, California 92521\\
USA\\
\printead{e2}}
\address[C]{R.~Y. Liu\\
Department of Statistics\\
Rutgers University\\
Piscataway, New Jersey 08854\\
USA\\
\printead{e3}}
\end{aug}

%
\received{\smonth{9} \syear{2014}}
%
\revised{\smonth{6} \syear{2015}}

%
\begin{abstract}
Statistical depth measures the
centrality of a point with respect to a given distribution or data
cloud. It provides a natural center-outward ordering of multivariate
data points and yields a systematic nonparametric multivariate analysis scheme.
In particular, the half-space depth is shown to have many desirable
properties and broad applicability.
However, the empirical half-space depth is zero outside the convex hull
of the data. This property has rendered the empirical half-space depth
useless outside the data cloud, and limited its utility in applications
where the extreme outlying probability mass is the focal point, such as
in classification problems and control charts with very small false
alarm rates. To address this issue, we apply extreme value statistics
to refine the
empirical half-space depth in ``the tail.''
This provides an important linkage between data depth, which is useful
for inference on centrality, and extreme value statistics, which is
useful for inference on extremity. The refined empirical half-space
depth can thus extend all its utilities beyond the data cloud, and
hence broaden greatly its applicability.
The refined estimator is shown to have substantially improved upon the
empirical estimator in theory and simulations. The benefit of this
improvement is also demonstrated through the applications in
classification and statistical process control.
\end{abstract}

%
\end{frontmatter}

\section{Introduction}\label{sec1}

Statistical depth generally is a measure of centrality with respect to
a multivariate distribution or a data cloud. It is shown to have many
useful data-driven features for developing statistical inference
methods and applications. For example, among other features, it can
also yield a center-outward ordering, and thus order statistics and
ranks for multivariate data.
With its rapid and broad advances, statistical depth has emerged to be a
powerful alternative approach in multivariate analysis.

There exist many different notions of statistical depth; see, for
example, \citet{LiuParSin99} and \citet{ZuoSer00} and
the references therein.
But the so-called geometric depths such as the half-space depth \citet{Tuk75} and the simplicial depth \citet{Liu90} are often preferred in many
nonparametric inference methods and applications for their intrinsic
desirable properties, as seen in \citet{DonGas92}, Liu and Singh
(\citeyear{LiuSin93,LiuSin97}), \citet{YehSin97}, \citet{RouHub99}, \citet{LiuParSin99}, \citet{ZuoSer00}, \citet{LiLiu04},
\citet{HalPaiSim10} and many others.

In practice, the empirical versions of the half-space depth and the
simplicial depth, however, suffer from the problem of vanishing value
outside the convex hull of the data. This problem is inherent in any
depth function that uses empirical counts based on the data to compute
its value. It renders the empirical version of such a depth useless
outside the data cloud, and limits its utility in applications
involving extreme outlying probability mass. A successful resolution to
this problem can avert such limitations and greatly enhance the utility
of depth functions. In investigating this problem, we observe that the
half-space depth involves projecting data points onto unit vectors, and
thus naturally lends itself in the framework of extreme value theory.
Therefore, we propose to refine the empirical half-space depth by
applying extreme value statistics to ``the tail.'' The aim of this paper
is to present this proposal, and assess and demonstrate the improvement
achieved by the proposal, in theory and applications.


To be more precise, let $\X_1, \ldots, \X_n$ be i.i.d. random vectors
taking values
in $\mathbb{R}^d,   d\geq1$. Denote the common probability
measure with $P$ and the empirical measure with $P_n$; denote
closed half-spaces with $H$. Then the half-space depth at $\x\in
\mathbb{R}^d$ is defined by
\[
D(\x)=\inf_{H:  \x\in H} P(H).
\]
Observe that the infimum can be restricted to half-spaces $H$
with $\x$ on their boundary. We can also write
\[
D(\x)=\inf_{\Vert \u\Vert =1}\mathbb{P} \bigl(\u^T
\X_1\geq\u^T \x\bigr),
\]
with $\Vert \cdot\Vert $ the radius or $L_2$-norm of a vector. The
classical nonparametric way to estimate $D(\x)$ is with the empirical
half-space depth:
\[
D_n(\x)=\inf_{H  :  \x\in H} P_n(H)=
\frac{1}{n} \inf_{\Vert \u\Vert =1} \# \bigl\{i\in\{1, \ldots, n\} :
\u^T \X_i\geq\u^T \x\bigr\}.
\]

It follows that for any $\x$ outside the convex hull of the data
$D_n(\x)=0$. This might seem a minor problem. Indeed, when the data are
univariate, the probability that a new observation falls outside the
convex hull is at most $2/(n+1)$, but in higher dimensions this
probability can be quite sizable. For example, for the multivariate
normal distribution and $n=100$ this probability is 8.8\% in dimension
2 and 21.7\% in dimension 3. Even when $n$ is as large as 500, these
probabilities are still 2.1\% ($d=2$) and 6.5\% ($d=3$); see, for
example, \citet{Efr65}. Outside the data hull, $D_n$ makes no
distinction between different points and provides hardly information
about $P$. This inability of distinguishing points in a sizable
subspace can severely restrict the utility of half-space depth in many
of its applications, such as statistical process control and
classification (see Section~\ref{sec3}). Note that the problem is not restricted
to $D_n$ being exactly 0: if $D_n(\x)$ is positive but very small, it
might not
adequately estimate $D(\x)$ due to the scarcity of useful data points.
Somewhat related, due to the discrete nature of $D_n$, ties occur
often. For example, $D_n(\X_i)=1/n$ for \textit{all} the data on the boundary
of the data hull, that is, all these data form one tie and cannot be
ranked effectively. (For the normal distribution in dimension 3 and
$n=500$ this tie, on average, has a size of about 32.) This phenomenon
renders rank procedures based on depth less precise and less efficient.

The goal of this paper is to refine the definition of empirical
half-space depth $D_n$ in the tail, that is, for values $\x$
where $D_n(\x)$ is zero or quite small. The proposed refined estimator
will be called $R_n$ (see Section~\ref{sec2} for the definition) and is based on
extreme value theory. The estimator $R_n$ is equal to $D_n$ in the
central region, where the depth is relatively high. Outside this region
$R_n$ is positive, smooth and it improves substantially on $D_n$.
Therefore, the aforementioned weaknesses of $D_n$ are ``repaired.''

As an illustration, we consider the estimation of the depth contour at
level $1/n$, that is, we want to estimate the set $\{\x\in\mathbb
{R}^d: D(\x)=1/n\}$, based on a random sample of size $n$. Using $D_n$,
it is usually estimated with the boundary of the data hull, where
indeed $D_n=1/n$, almost surely. We also estimate it using our refined
estimator by $\{\x\in\mathbb{R}^d: R_n(\x)=1/n\}$.
We consider as an example the bivariate spherical Cauchy distribution
(see Section \ref{sec2.3}) and simulate one random sample of size $n=500$; see
Figure~\ref{bestplotever}. (The computation of these depth contours is
discussed in Remark~\ref{re6} of Section~\ref{sec2.2}.) It clearly shows that $R_n$
greatly improves $D_n$; $D_n$ fails completely here, whereas $R_n$
performs well. This indicates that our refined estimator can be very
useful in practice.

%

In the next section, we will define $R_n$ and show, under appropriate
conditions, its uniform ratio consistency (considering $R_n/D-1$) on a
very large region, much larger than the data hull. In contrast, $D_n/D$
is not uniformly close to 1 \textit{on} the data hull. We further show
through simulations that these asymptotic differences between $R_n$ and
$D_n$ are clearly present for finite samples, that is, that $R_n$
substantially outperforms $D_n$ in the tail. In Section~\ref{sec3}, we
investigate the impact of these theoretical improvements in real
applications of data depth using examples in statistical process
control (SPC) and classification. Both applications obtain substantial
improvements by using $R_n$.
Finally, we provide some concluding remarks in Section~\ref{sec4}.
All proofs are deferred to Section~\ref{sec5}.

\section{Methodology and main results}\label{sec2}
\subsection{Dimension one}\label{sec2.1}

We first consider refining $D_n$ in the one-dimensional case,
particularly since it serves as a building block for us to refine $D_n$
in higher dimensions. Let $X_1, \ldots, X_n$ be i.i.d. random
variables with common continuous distribution function $F$ with $0<F(0)<1$.
Write $S=1-F$. Let $F_n$ be the (right-continuous) empirical
distribution function and define $S_n(x)=1-F_n(x^-)$. The half-space
depth and its empirical counterpart in the one-dimensional case are
simply $D(x)=\min(F(x), S(x))$ and $D_n(x)=\min(F_n(x), S_n(x))$,
respectively. It is clear that the aforementioned shortcomings of $D_n$
are due to the inadequacy of the empirical distribution function as an
estimator in the tails. Since extreme value statistics is well suited
for inference problems in this setting, we propose applying it to
refine $D_n$ in the tails.

%
\begin{figure}

\includegraphics{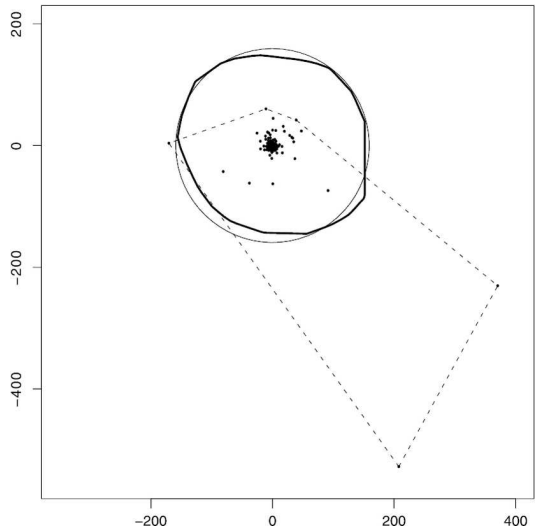}

\caption{Depth contours at level $1/n$ based on $D$ (circle), $D_n$
(dashed) and $R_n$ (solid) for a standard bivariate spherical Cauchy
random sample; $n=500$.}\label{bestplotever}
\end{figure}

In extreme value theory, it is assumed that there exist a location
function $b$ and a scale function $a>0$ such that
%
\begin{eqnarray}
\label{aaa} \lim_{t \to\infty} t \bigl(1 - F\bigl(a(t) y + b(t)
\bigr) \bigr) = - \log G_\gamma (y)= (1+\gamma y)^{-1/\gamma},
\nonumber
\\[-8pt]
\\[-8pt]
\eqntext{1+\gamma
y > 0.}
\end{eqnarray}
Here, $G_{\gamma}$ is the limiting extreme value distribution and
$\gamma\in\mathbb{R}$ is the extreme value index. If (\ref{aaa})
holds, $F$ is said to be in the max domain of attraction of $G_{\gamma
}$. See, for example, \citet{deHFer06}. The above assumption
guarantees that $F$ has a ``regular'' tail and makes extrapolation
outside the data range possible.

If $F$ is in the max-domain of attraction of $G_{\gamma}$, by setting
$t=n/k$ and $x = a(t) y + b(t)$ in (\ref{aaa}), we obtain for large
$n/k$ and large $x$
%
\begin{equation}
\label{eq2_4} \mathbb{P}(X >x)=1-F(x) \approx\frac{k}{n} \biggl(1+\gamma
\frac{x- b(n/k)}{a(n/k)} \biggr)^{-1/\gamma}.
\end{equation}
%
Let $\hat\gamma$ and $\hat a=\hat a(n/k)$ be estimators for
$\gamma$ and $a=a(n/k)$, respectively. Define $\hat b=\hat
b(n/k)=X_{n-k:n}$, where $X_{i:n}$ denotes the $i$th order
statistic of $X_1, \ldots, X_n$.
Plugging these estimators into (\ref{eq2_4}), we obtain the following
estimator for the right-tail probability $1-F(x)$:
%
\begin{equation}
\label{eq4} p^r_n(x) = \frac{k}{n} \biggl(\max
\biggl[0,1+\hat{\gamma}\frac{x- \hat{b}}{\hat
{a}} \biggr] \biggr)^{-1/\hat{\gamma}}.
\end{equation}
To estimate the left-tail probability, we can define $p^l_n(x)$
similarly as $p^r_n(x)$ by using the $-X_i$.

The general idea of estimating $D$ with our refined estimator is the
following. For a given $k$, we define the central region to be
$(X_{k+1:n},X_{n-k:n})$. For $x$ in this central region, we define
$R_n(x)=D_n(x)$, that is, we use the classical empirical half-space
depth. In the right tail, that is, when $x\geq X_{n-k:n}$, we refine
$D_n$ by defining $R_n(x)=p^r_n(x)$ and similarly, when $x\leq
X_{k+1:n}$ (the left tail), we set $R_n(x)=p^l_n(x)$.
At the ``glue-up'' points $X_{n-k:n}$ and $X_{k+1:n}$, we have
\begin{eqnarray*}
R_n(X_{n-k:n})&=&p^r_n(X_{n-k:n})=
\frac
{k}{n}=1-F_n(X_{n-k:n})=D_n
\bigl(X_{n-k:n}^+\bigr),
\\
R_n(X_{k+1:n})&=&p^l_n(X_{k+1:n})=
\frac{k}{n}=F_n\bigl(X_{k+1:n}^-\bigr)=D_n
\bigl(X_{k+1:n}^-\bigr).
\end{eqnarray*}

In the following, we study the asymptotic properties of our refined
empirical half-space depth $R_n$. Throughout we assume that $k=k_n<n/2$
is an intermediate sequence: a sequence of positive integers
satisfying
%
\begin{equation}
\label{bbb}k\to\infty \quad\mbox{and}\quad k/n\to0 \qquad\mbox{as } n \to\infty.
\end{equation}
We need a second-order condition in both the left tail and the right
tail; for simplicity, we will only specify it for the right tail. Let
$V(t)=F^{-1}(1-1/t),  t>1$, be the tail quantile function. We
can and will take the location function $b(t)=V(t)$. We assume
that the derivative $V'$ exists and that for some eventually
positive or eventually negative function $A$ with
$\lim_{t\to\infty}A(t)=0$ and for some $\rho< 0$ we have
%
\begin{equation}
\label{ccc}\lim_{t\to\infty} \frac{
{V'(tx)}/{V'(t)}-x^{\gamma-1}}{A(t)}=x^{\gamma-1}
\frac{x^{\rho}-1}{\rho
}, \qquad x>0.
\end{equation}
This condition implies (for $\rho<0$)
%
\begin{equation}
\label{finite} \lim_{t\to\infty}\sup_{y\geq1/2,  y\neq1} \biggl
\llvert \frac{
({(V(ty)-V(t))}/{(tV'(t))})({\gamma}/{(y^\gamma-1)})-1}{A(t)}\biggr\rrvert <\infty.
\end{equation}
This limit relation is somewhat similar to Lemma~4.3.5 in \citet{deHFer06}. A proof can be given along the lines
of the proof of that lemma; the proof uses in particular
Theorem~2.3.9 in \citet{deHFer06}, with $U$ and
$\gamma$ there replaced by $V'$ and $\gamma-1$, respectively.
We can and will take the scale function $a(t)=tV'(t)$. We
assume
%
\begin{equation}
\label{ddd}\sqrt{k}A(n/k)\to\lambda\qquad \mbox{for some } \lambda\in\mathbb{R}.
\end{equation}
We will also assume that the estimators $\hat\gamma$ and $\hat
a$ are such that
%
\begin{equation}
\label{eee} \Gamma_n:=\sqrt{k}(\hat\gamma-\gamma)=O_p(1)\quad
\mbox{and} \quad \sqrt{k} \biggl(\frac{\hat
a}{a}-1 \biggr)=O_p(1).
\end{equation}
This condition is known to hold for various estimators of $\gamma$ and $a$;
see de Haan and Ferreira [(\citeyear{deHFer06}), Chapters~3 and 4]. Define
\[
w_\gamma(t)=t^{-\gamma}\int_1^t
s^{\gamma-1}\log s \,ds,\qquad  t>1.
\]
Note that, as $t\to\infty$,
\[
w_\gamma(t)\sim\cases{
\displaystyle\frac{1}{\gamma}
\log t, & \quad $\gamma>0$,
\vspace*{2pt}\cr
\displaystyle\frac{1}{2} ( \log t ) ^{2}, &\quad $\gamma=0$,
\vspace*{2pt}\cr
\displaystyle\frac{1}{\gamma^{2}}t^{-\gamma}, &\quad  $\gamma<0$.}
\]

\begin{theorem}
\label{th:1d}
Let $\delta_n$ be a sequence of numbers in $(0,1/2)$ such that
$n\delta_n\to0$ as $n\to\infty$. Assume that (\ref{aaa}) and its
left-tail counterpart hold; also assume
$w_\gamma (\frac{k}{n\delta_n} )/\sqrt{k}\to0$ as
$n\to\infty$. Then, if (\ref{bbb}), (\ref{ccc}), (\ref{ddd}) and (\ref
{eee}) hold, we have
\[
\sup_{x\in\mathbb{R}:  D(x)\geq\delta_n } \biggl\llvert \frac{R_n(x)}{D(x)}-1\biggr\rrvert
\stackrel{p} {\to} 0\qquad \mbox{as } n \to\infty.
\]
\end{theorem}

The condition on $\delta_n$ and $k$ specializes to
\[
\frac{\log(n \delta_n)}{\sqrt k}\to0\qquad \mbox{for } \gamma>0\quad \mbox{and}\quad \frac{\log^2(n \delta_n)}{\sqrt k}\to0
\qquad\mbox{for } \gamma=0.
\]

\begin{rem}\label{re1} The main focus of this paper is on the
tails where both $R_n$ and $D$ are small, and as such, $R_n-D$ (just
like $D_n-D$) is inherently small as well. This implies that the usual
consistency statement $\sup_x |R_n(x)-D(x)|\stackrel{p}{\to}0$ is not
particularly meaningful for assessing the performance of $R_n$ as an
estimator of $D$. Instead, we consider the ratio consistency in terms
of $R_n/D-1$ as stated in Theorem~\ref{th:1d}.
Note that, in addition to the usual consistency $\sup_x
|D_n(x)-D(x)|\stackrel{p}{\to}0$ [\citet{DonGas92}],
Theorem~\ref{th:1d}
also holds for $D_n$, when $n\delta_n\to\infty$, but
\textit{not} when $n\delta_n$ tends to a nonnegative constant; cf.
(\ref{sw}) below.
This shows that the region for which $R_n/D$ is close to 1 (for large
$n$ and with high probability) is much greater than that for $D_n/D$.
\end{rem}

\begin{rem}\label{re2}
It is natural to consider an asymptotic normality result instead of
the consistency result in Theorem~\ref{th:1d}, but note that the convergence rate
($1/r_n$, say, with, $r_n/\sqrt{n}\to0$) for the process $R_n/D-1$ in
such a result will be determined by $x_n$-values with $D(x_n)\to0$; at
a fixed $x$ the weak limit of
$r_n(R_n(x)/D(x)-1)=( r_n/\sqrt{n})\sqrt{n}(R_n(x)/D(x)-1)$ will be 0.
This means that a proper refinement of Theorem~\ref{th:1d}, specifying the rate
of convergence and providing a nondegenerate limit, is not possible. On
the other hand, if we consider a single $x=x_n$ in the right tail such
that $nD(x_n)/k\to0$, then it follows from Theorem~4.4.1 in \citet{deHFer06} (under the assumptions there) that for some $\mu$
and $\sigma>0$
\[
\frac{\sqrt{k}}{w_\gamma(k/(nD(x_n)))} \biggl(\frac
{R_n(x_n)}{D(x_n)}-1 \biggr)\stackrel{d} {\to} N\bigl(
\mu, \sigma^2\bigr) \qquad\mbox{as } n\to\infty,
\]
since $R_n(x_n)=p_n^r(x_n)$, see (\ref{eq4}), with probability tending
to one. Indeed, the convergence rate here is slower than for fixed $x$:
$\frac{r_n}{\sqrt{n}}=\frac{\sqrt{k}}{w_\gamma(k/(nD(x_n)))\sqrt n}\to0$.
\end{rem}

\subsection{Higher dimensions}\label{sec2.2}

We next consider constructing the refined half-space depth
estimator in the more interesting, multivariate case, that is, $d\geq2$.
Let $\X_1, \ldots, \X_n$ be i.i.d. random vectors drawn from a common
continuous distribution function $F$.
To refine $D_n$ we need now some more structure for $F$. More
precisely, we assume multivariate
regular variation for $F$, that is, there exists a measure $\nu$ such that
%
\begin{equation}
\label{nu} \lim_{t\rightarrow\oo}\frac{\mathbb{P} ( \X_1 \in tB)}{\mathbb{P} (\Vert \X
_1\Vert  \geq t)}=\nu(B)<\infty,
\end{equation}
for every Borel set $B$ on $\mathbb{R}^d$ that is bounded away
from the origin and satisfies $\nu(\partial B)=0$; see, for example,
\citet{JesMik06}. Note that the choice of the ``spherical''
$L_2$-norm is not relevant: any other norm can be used instead.
This implies that for some $\alpha>0$
\[
\label{rad} \lim_{t\rightarrow\oo}\frac{\mathbb{P} ( \Vert \X_1\Vert \geq tx )}{\mathbb{P}
(\Vert \X_1\Vert  \geq t)}=x^{-\alpha}\qquad
\mbox{for } x>0.
\]
The parameter $\alpha$ is called the tail index and $\gamma=1/\alpha>0$
is the extreme value index. Note that, for all $a>0$, $\nu
(aB)=a^{-\alpha} \nu(B)$. We further require
%
\begin{equation}
\label{radi}\frac{\mathbb{P}(\Vert \X_1\Vert >t)}{t^{-\alpha
}}\to c\in(0,\infty).
\end{equation}
This simple condition in effect replaces the second-order condition of
the univariate case, although it is a slightly weaker condition; cf.
Cai, Einmahl and de~Haan (\citeyear{CaiEindeH11}), page 1807.
We also assume that
%
\begin{equation}\quad
\label{fff} \u^T \X_1 \mbox{ has a continuous
distribution function } F_\u\mbox{ for every unit vector } \u,
\end{equation}
and that, with $H_{r,\u}:=\{\x\in\mathbb{R}^d  :   \u^T \x\geq r \}
,  r>0$,
%
\begin{equation}
\label{inf}\inf_{\Vert \u\Vert =1} \nu(H_{1,\u}) >0.
\end{equation}
Note that the continuity of the $F_\u$ implies the continuity of $D$.
Also, observe that
the multivariate regular variation condition (\ref{nu}) implies that
for every unit vector $\u$, $F_\u$ is in the univariate max domain of
attraction with the same $\gamma=1/\alpha$:
as $t\to\infty$,
\begin{eqnarray*}
\frac{1-F_\u({tr})}{1-F_\u(t)}&=&\frac{\mathbb{P} ( \X_1 \in trH_{1,\u
})}{\mathbb{P} (\X_1 \in tH_{1,\u})}= \frac{\mathbb{P} ( \X_1 \in trH_{1,\u})}{\mathbb{P} (\Vert \X_1\Vert  \geq
t)}\cdot\frac{\mathbb{P} (\Vert \X_1\Vert  \geq t)}{\mathbb{P} (\X_1 \in
tH_{1,\u})}
\to\frac{\nu(rH_{1,\u})}{\nu(H_{1,\u})}\\
&=&r^{-\alpha}.
\end{eqnarray*}

Recall that the half-space depth, relative to $\mathbb{P}$, is defined as
\[
D(\x)=\inf_{\Vert \u\Vert =1}\mathbb{P} \bigl(\u^T
\X_1\geq\u^T \x\bigr).
\]
To estimate $D(\x)$, we only need to estimate the one-dimensional tail
probabilities $\mathbb{P} (\u^T \X_1\geq\u^T \x)$ along each
projection direction $\u$. Since we already know how to construct the
refined estimator for a tail probability in the one-dimensional case,
we are now ready to define our refined empirical half-space depth $R_n$
in dimension~$d$.

More specifically, fix a direction (a unit vector) $\u$. Consider the
univariate data $W_i=\u^T \X_i$, $i=1, \ldots, n$. We can refine the
tail probability estimator of the $W_i$ similarly as in the previous
subsection, but since $\gamma>0$ we can use $a =\gamma b$. This leads,
for $w\geq W_{n-k:n}$,
to a simplified estimator of the right-tail probabilities:
%
\begin{equation}
\label{eq:md} p_{n,\u}(w)=\frac{k}{n} \biggl(\frac{w}{W_{n-k:n}}
\biggr)^{-\widehat
\alpha};
\end{equation}
cf. (\ref{eq2_4}) and (\ref{eq4}).
The estimator $\widehat\alpha=1/\hat\gamma$ will
be based on the $\Vert \X_i\Vert $. We will assume that
$\hat\gamma$ is such that
%
\begin{equation}
\label{gamma} \Gamma_n:=\sqrt{k}(\hat\gamma-
\gamma)=O_p(1).
\end{equation}
For $w < W_{n-k:n}$
an estimator of $1-F_\u(w)$ is simply $1-F_{n,\u}(w)$, with $F_{n,\u}$
the empirical distribution function of $W_1,\ldots,W_n$.
Denote the thus obtained estimator of $1-F_\u$ with $
1-\widehat F_{\u}$.
This leads
to the refined estimator of $D(\x)$:
\[
R_n(\x)= \inf_{\Vert \u\Vert =1} 1-\widehat F_{\u}
\bigl(\u^T\x-\bigr).
\]
%

Next, we present the analogue of Theorem~\ref{th:1d} for the
multivariate $R_n$. Note that it is much more complicated to analyze
$R_n$ here than in dimension one, since for every $\x\in
\mathbb{R}^d$ we have infinitely many directions $\u$ instead of
only two.

\begin{theorem}\label{2} 
Let $\delta_n$ be a sequence of numbers in $(0,1/2)$ such
that $n\delta_n\to0$ as $n\to\infty$. Also assume
$\log(n\delta_n)/\sqrt{k}\to0$ as $n\to\infty$. Then, if (\ref{nu}),
(\ref{bbb}), (\ref{radi}), (\ref{fff}), (\ref{inf}) and (\ref{gamma})
hold, we have
%
\begin{equation}
\label{theo2}\sup_{\x\in\mathbb{R}^d:  D(\x)\geq
\delta_n } \biggl\llvert \frac{R_n(\x)}{D(\x)}-1
\biggr\rrvert \stackrel{p} {\to} 0\qquad \mbox{as } n \to\infty.
\end{equation}
\end{theorem}

\begin{rem}\label{re3} It is known that the half-space depth is
affine invariant. This means that the depth value does not change under
a linear transformation. Specifically, $D(\x)=D_{\A,\b}(\A\x+\b)$,
where $D_{\A,\b}$ indicates the depth value based on the sample $\A\X_i
+\b, i=1, \ldots, n$. Here, $\A$ is a $d \times d$ nonsingular matrix
and $\b\in\mathbb{R}^d$. Although this property does not hold for $R_n$
exactly, it holds approximately through (\ref{theo2}).
\end{rem}

\begin{rem}\label{re4} The class of multivariate regularly varying
distributions [see (\ref{nu})] is quite broad. It contains, for
example, all elliptical distributions with a heavy tailed radial
distribution (such as multivariate $t$-distributions) and all
distributions in the sum domain of attraction of a multivariate
(nonnormal) stable distribution; see, for example, \citet{MeeSch01}, part III. Some examples are seen in Section~\ref{sec2.3}.
Note in particular that the extreme density contours of such
distributions can have more or less arbitrary shapes, not only spheres
or ellipsoids. Two such distributions, with nonconvex or asymmetric
extreme density contours, can be found in Cai, Einmahl and de~Haan (\citeyear{CaiEindeH11}). It is
also worth noting that the multivariate regular variation condition can
be verified using the test in \citet{EinKra15}.
\end{rem}

\begin{rem}\label{re5}
For $ D_n$ the statement
of Theorem~\ref{2} holds when $n\delta_n\to\infty$ [see (\ref{i})
below] but
\textit{not} when $n\delta_n$ tends to a nonnegative constant, which
again shows
that $R_n/D$ is close to 1 (for large $n$ and with high probability)
on a much larger region than where $D_n/D$ is.
\end{rem}

\begin{rem}\label{re6}
(i) \textit{Computation of $R_n$}: Recall
that when $D_n$ or $ R_n$ is at least $k/n$, then they are equal. Let
$\x$ be such that $D_n(\x)= R_n(\x)=k/n$ and let $\x^*=c\x$ with $c>1$.
Based on (\ref{eq:md}), we obtain
%
\begin{equation}
\label{eq:rm6} R_n\bigl(\x^*\bigr)= c^{-\widehat\alpha}R_n(
\x).
\end{equation}
Combination of both properties enables us to calculate $R_n$ readily by
utilizing any available algorithm for computing $D_n$.


(ii) \textit{Computation of depth contour based
on $R_n$ in Figure~\ref{bestplotever}}: Write $\x^*=(r_{\theta}\cos\theta,r_{\theta}\sin\theta
)$. We need to find $r_{\theta}$ such that $R_n(\x^*)=1/n$ for all
$\theta\in[0,2\pi)$. For any fixed $\theta$, similar to the above
procedure for computing $R_n$, we first find $\x=(s_{\theta}\cos\theta
,s_{\theta}\sin\theta)$ such that $D_n(\x)=k/n$. Then based on (\ref
{eq:rm6}), $\x^*=k^{1/\widehat\alpha}(s_{\theta}\cos\theta,s_{\theta
}\sin\theta)$ and $R_n(\x^*)=1/n$. The $R_n$-depth contour in Figure~\ref{bestplotever} is drawn using 500 evenly distributed $\theta$'s in
$[0,2\pi)$.
\end{rem}

\begin{rem}\label{re7} Our estimator $R_n$ involves $k$ and its
performance obviously will be affected by the choice of $k$. The
problem of choosing optimal $k$ is an inherent one in extreme value
statistics. Various approaches have been proposed in the literature.
One commonly used heuristic approach is to plot the relevant estimator
versus $k$, visually identify the first (or earliest) stable
(approximately constant) region in the plot, and then choose the
midpoint of this region as $k$.
This approach is the one we used in our numerical studies below. We
find more or less the same value of $k$ in the first few samples. For
some specific problem settings, procedures for determining the optimal
$k$ have been developed. It would be worthwhile developing such a
procedure for $R_n$. Meanwhile, we note that even with the present
choice of $k$, which may well be only suboptimal, $R_n$ already clearly
outperforms~$D_n$.
\end{rem}

\subsection{Simulation comparison between $R_n$ and $D_n$}\label{sec2.3}

In this section, we\break present a simulation study to compare the
performance of our refined empirical half-space depth $R_n$ with the
performance of the original empirical half-space depth $D_n$. We
consider the following distributions in our simulation study:\vadjust{\goodbreak}
\begin{itemize}
\item Standard normal distribution. This is a light-tailed
distribution with $\gamma=\rho=0$.
\item Cauchy distribution. This is a very heavy-tailed
distribution with $\gamma=1$ and $\rho=-2$.\label{page11}
\item $t$-distribution with 2 degrees of freedom. This is a
heavy-tailed distribution with $\gamma=1/2$ and $\rho=-1$.
\item Burr-type distribution, which is a symmetric
distribution about 0 with density
\[
f(x)=\frac{3|x|^5}{2(1+x^6)^{3/2}} , \qquad x \in\mathbb{R}.
\]
This distribution is less heavy-tailed with $\gamma=1/3$ and $\rho=-2$.
\item Standard bivariate normal distribution. This is a
light-tailed distribution with $\gamma=0$.
\item Bivariate spherical Cauchy distribution with density
\[
f(x,y)=\frac{1}{2\pi(1+x^2+y^2)^{3/2}} , \qquad (x,y) \in\mathbb{R}^2.
\]
This is a very heavy-tailed distribution with $\gamma=1$.
\item Bivariate elliptical distribution with density ($r_0
\approx1.2481$)
\[
f(x,y)= %
\cases{ \displaystyle\frac{3}{4\pi}r_0^4
\bigl(1+r_0^6\bigr)^{-3/2}, &\quad
$x^2/4+y^2<r_0^2$,\vspace*{2pt}
\cr
\displaystyle\frac{3(x^2/4+y^2)^2}{4\pi(1+(x^2/4+y^2)^3)^{3/2}} , &\quad  $x^2/4+y^2 \geq
r_0^2$.} %
\]
This is a less heavy-tailed distribution with $\gamma=1/3$.
\item
Bivariate ``clover'' distribution with density ($r_0\approx1.2481$)
%
\begin{equation}
\label{clov} f(x,y)=\cases{ \displaystyle
\frac{3}{10\pi}r_0^4
\bigl(1+r_0^6\bigr)^{-3/2} \biggl(5+
\frac
{4(x^2+y^2)^2-
32x^2y^2}{r_0(x^2+y^2)^{3/2}} \biggr),\vspace*{2pt}\cr
\qquad x^2+y^2< r_0^2,
\vspace*{2pt}\cr
\displaystyle\frac{3 (9(x^2+y^2)^2-32x^2y^2 )}{10\pi
(1+({x^2}+y^2)^3)^{3/2}} ,\vspace*{2pt}\cr
\qquad x^2+y^2\geq r_0^2.}
\end{equation}
This is again a less heavy-tailed distribution with $\gamma=1/3$,
however, it is not an elliptical distribution; see Cai, Einmahl and de~Haan (\citeyear{CaiEindeH11}).
\item Trivariate spherical Cauchy distribution with density
\[
f(x,y,z)=\frac{1}{\pi^2(1+x^2+y^2+z^2)^2} ,\qquad (x,y,z) \in\mathbb{R}^3.
\]
This is a very heavy-tailed distribution with $\gamma=1$.
\item Quadrivariate spherical Cauchy distribution with density
\[
f(x,y,z,w)=\frac{3}{4\pi^2(1+x^2+y^2+z^2+w^2)^{3/2}} ,\qquad (x,y,z,w) \in\mathbb{R}^4.
\]
This is again a very heavy-tailed distribution with $\gamma=1$.
\end{itemize}
The first four distributions will be used to assess the finite sample
performance of Theorem~\ref{th:1d} (although for the standard normal
distribution $\rho<0$ does not hold), and the last five distributions
are used to assess the finite sample performance of Theorem~\ref{2}.

For each of the above distributions, we first generate a random sample
of size 500. Based on this random sample, $R_n$ and $D_n$ are then
calculated for a point $\x$ where the theoretical depth $D(\x)$ is
$1/100$, $1/500$, $1/1000$ and $1/2000$, respectively.
To calculate $R_n$, an estimator of $\gamma=1/\alpha$ (and $a$) is
needed. For the univariate distributions, we use the moment estimator
of \citet{DekEindeH89} for estimating $\gamma$ and for
$a$ we use a corresponding estimator; see formula (4.2.4) in \citet{deHFer06}. For the multivariate distributions (except the
bivariate normal), since we assume that $\gamma>0$, we use the \citet{Hil75} estimator, based on the $\Vert \X_i\Vert $. For the bivariate normal
distribution, because it does not satisfy the conditions of Theorem~\ref
{2}, we use (\ref{eq4}) instead of (\ref{eq:md}) to estimate the
right-tail probability of the $W_i$. In other words,
\[
p_{n,\u}(w)=\frac{k}{n} \biggl(\max \biggl[0,1+\hat\gamma
\frac{w-\hat
{b}_{\u}}{\hat{a}_{\u}} \biggr] \biggr)^{-1/\hat\gamma},
\]
where $\hat{\gamma}$ is the moment estimator based on the $\Vert \X_i\Vert $,
$\hat{b}_{\u} = W_{n-k:n}$, and $\hat{a}_{\u}$ is again as in
(4.2.4) in \citet{deHFer06}. Since (\ref{eq:rm6}) does not
hold for this case any more, we follow \citet{CueNie08} to approximate $R_n$ using 500 $\u$'s that are uniformly and
independently distributed on the unit sphere. For all 10 distributions
the value of $k$ is selected by searching visually for the first stable
part in the plots, based on 3 to 5 samples, as described in more detail
in Remark~\ref{re7}. This leads to values of $k$ ranging from 50 to 100: 6 times
50, twice 75 and twice 100.

\begin{figure}

\includegraphics{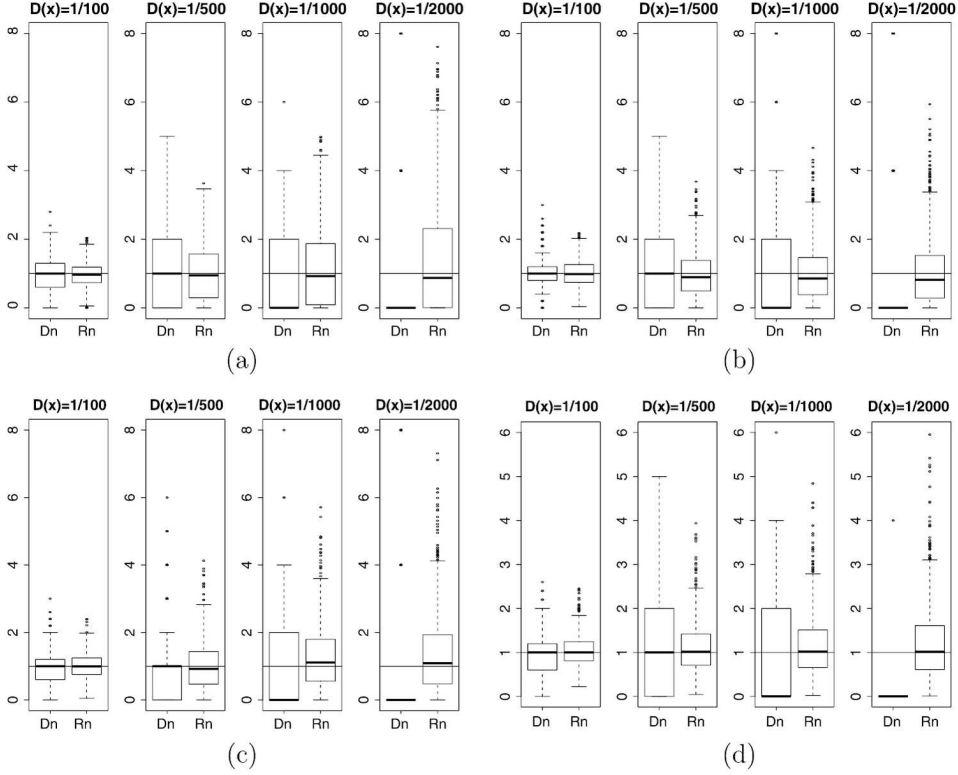}

\caption{Comparison of $D_n/D$ (left) and $R_n/D$ (right) at 4
decreasing levels under \textup{(a)} normal distribution; \textup{(b)} Burr-type
distribution; \textup{(c)} $t$-distribution with 2 degrees of freedom; \textup{(d)}
Cauchy distribution.}
\label{fig:1D}
\end{figure}

%
%

\begin{figure}

\includegraphics{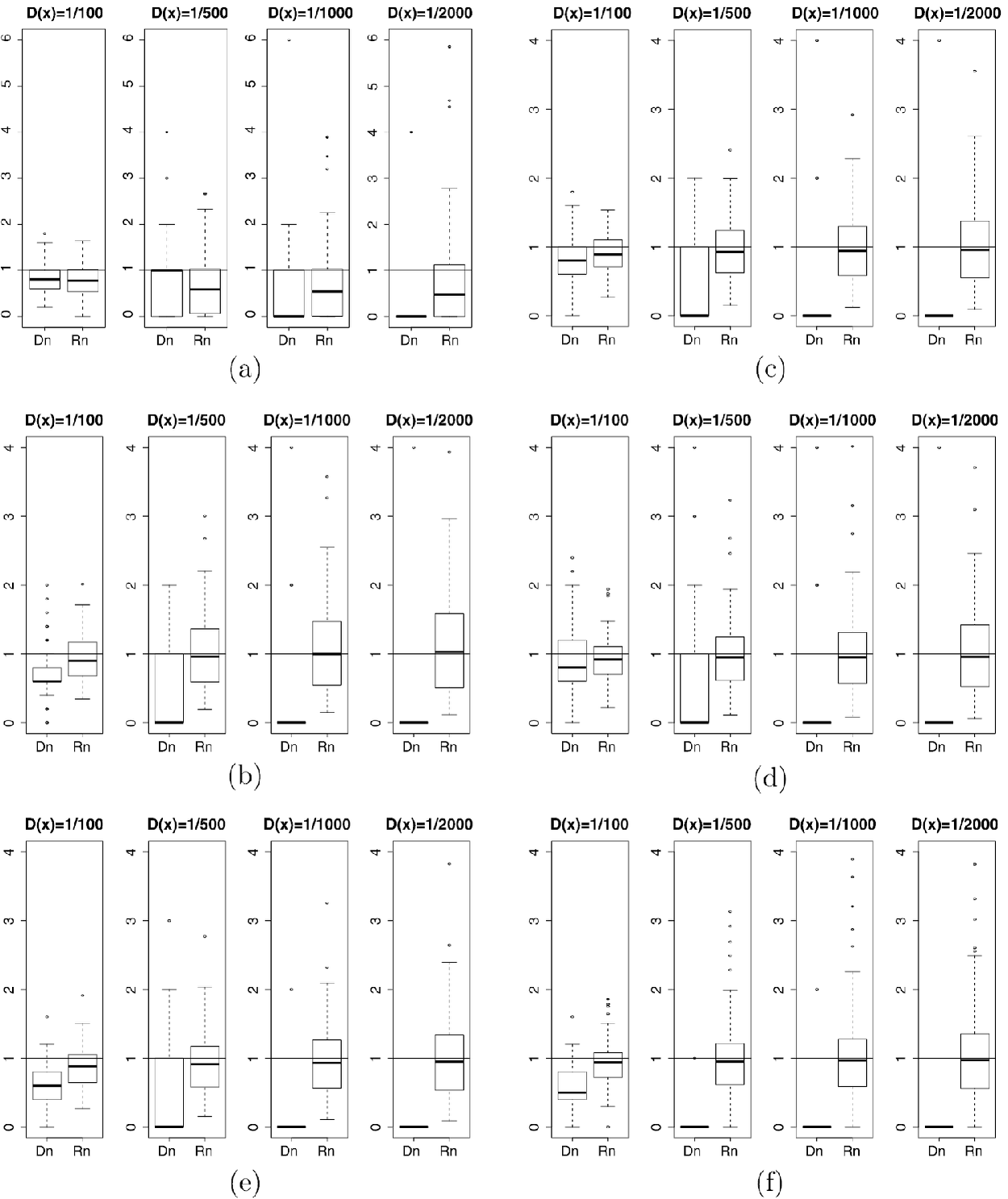}

\caption{Comparison of $D_n/D$ (left) and $R_n/D$ (right) at 4
decreasing levels under \textup{(a)} bivariate normal distribution; \textup{(b)}
bivariate spherical Cauchy distribution; \textup{(c)} bivariate elliptical
distribution; \textup{(d)} bivariate clover distribution; \textup{(e)} trivariate
spherical Cauchy distribution; \textup{(f)} quadrivariate spherical Cauchy distribution.}
\label{fig:2D}
\end{figure}

%

We carry out the above simulation 100 times for each of the
distributions. The boxplots of $R_n(\x)/D(\x)$ and $D_n(\x)/D(\x)$ for
each of the four depth levels from the 100 simulations for different
distributions are plotted in Figures~\ref{fig:1D} and \ref{fig:2D}. As
we can see from those boxplots, for all the four depth levels and all
the distributions except the bivariate normal distribution, the $R_n(\x
)/D(\x)$ are all well centered at 1. In contrast, the original
empirical halfspace depth $D_n$ can only provide a reasonable estimate
of $D$ when $D$ is not too small. When $D(\x)$ is small relative to
$n$, most of the $D_n(\x)$ are zero. These results support the
theoretical findings that $R_n$ is a better estimator than $D_n$ in the
tail. For the bivariate normal distribution, although it does not
satisfy the assumptions of Theorem~\ref{2}, the performance of $R_n$ is
still much better than the performance of $D_n$.

\section{Impact of the refinement of $D_n$ on applications}\label{sec3}

\subsection{Statistical process control}\label{sec3.1}

In this section, we present two applications where $R_n$ significantly
improves the performance of the depth based procedures over $D_n$. The
first one is statistical process control (SPC). SPC is the application
of statistical methods to the monitoring of a process outcome in order
to detect abnormal variations of the process from a specified
in-control distribution. It has many applications in manufacturing
processes. A typical setup for SPC is the following. There are $n$
i.i.d. historical (reference) data for the monitored process outcome,
denoted by $\X_1, \ldots, \X_n$ $ \in\mathbb{R}^d$ $ (d\geq1)$, from
the in-control process. Let $F_0$ be the underlying distribution of the
$\X_i$, also referred to as the in-control distribution. Let $\Y_1,\Y
_2, \ldots$ be future observations of the process outcome, under the
distribution $F_1$. The task of SPC is to determine if $F_1$ is the
same as $F_0$ and if not, to signal when $F_1$ changes from $F_0$ as
early as possible.


When the process outcome is multivariate and follows a multivariate
normal distribution, an SPC procedure with a false alarm rate $\alpha$
can be defined as follows: $\Y_i$ is out of control if
$T_i^2>{d(n+1)(n-1)/(n(n-d))}F_{d,n-d}(\alpha)$, where $T_i^2=(\Y_i
-\bar{\X})'S^{-1}(\Y_i -\bar{\X})$, $\bar{\X}=\sum_{i=1}^n \X_i/n$,
$S=\sum_{i=1}^n(\X_i-\bar{\X})(\X_i-\bar{\X})'/(n-1)$, and
$F_{d,n-d}(\alpha)$ is the upper $\alpha$ quantile of an $F$
distribution with $d$ and $n-d$ degrees of freedom.

The above procedure requires that the process outcome follows a
multivariate normal distribution. Therefore, we refer to it as the
parametric SPC procedure hereafter. In many real world applications,
the normality assumption may not hold. Therefore, a nonparametric SPC
procedure is more desirable. Following \citet{Liu95}, a nonparametric SPC
procedure with a false alarm rate $\alpha$ can be defined as follows:
$\Y_i$ is out of control if $\#\{\X_j  :   D(\Y_i)>D(\X_j), j=1,
\ldots, n\}/n < \alpha$, where $D$ is the depth with respect to $F_0$.
Since the in-control distribution is usually unknown in practice, $D$
in the above procedure is usually replaced by $D_n$, the empirical
depth with respect to the historical data, $\X_1, \ldots, \X_n$.

Due to its completely nonparametric nature and its capability of
characterizing the geometric structure of the underlying distribution,
the half-space depth is a popular choice in the above depth based SPC
procedure. Because the future process outcomes $\Y_i$ that lie in the
outskirts of the historical data are more of concern in this SPC
procedure, how close the achieved false alarm rate to the nominal level
$\alpha$ depends on how well the empirical half-space depth $D_n$
estimates the theoretical half-space depth $D$ for those points. As
shown in this paper, this estimation is not satisfactory when $n$ is
not large enough. Therefore, the achieved false alarm rate can severely
deviate from its nominal level $\alpha$ when
$D_n$ is used. To overcome this drawback of using $D_n$, we use our
refined halfspace depth $R_n$ in the above SPC procedure instead. Based
on the results in Section~\ref{sec2}, we expect the above depth based SPC
procedure will achieve the nominal false alarm rate if $R_n$ is used.

\begin{figure}

\includegraphics{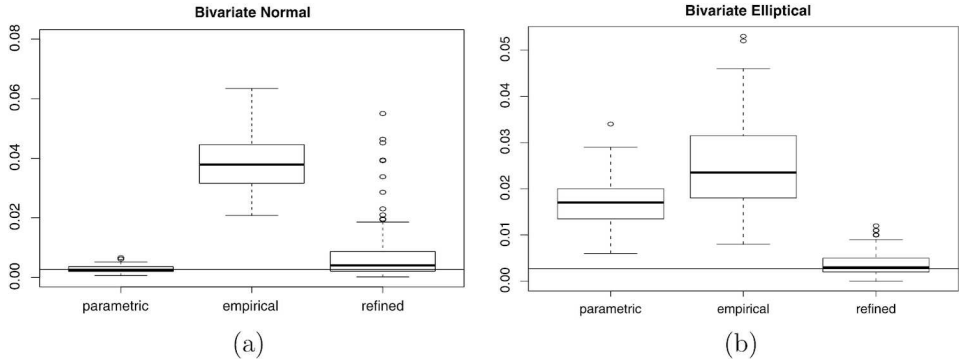}

\caption{The achieved false alarm rates for the parametric procedure,
the $D_n$ based procedure and the $R_n$ based procedure under \textup{(a)}
bivariate normal distribution; \textup{(b)} bivariate elliptical distribution.}
\label{fig:QC0}
\end{figure}

To demonstrate the performance of the $R_n$ based SPC procedure, we
carry out the following simulation. We first generate $n=500$
historical data $\X_i$ from the standard bivariate normal distribution.
We then generate another $5000$ future observations $\Y_i$ from the
same bivariate normal distribution. We apply to the 5000 $\Y_i$ the
following three SPC procedures: the parametric procedure, the $D_n$
based procedure and the $R_n$ based procedure.
We calculate $R_n$ for the bivariate normal distribution as described
in the previous section.
The nominal false alarm rate $\alpha$ for each procedure is set to be
at 0.0027 (the false alarm rate for the popular 3-sigma procedure in
the univariate normal setting). The achieved false alarm rate for each
procedure is then calculated as the proportion of $\Y_i$ being labeled
as out-of-control by its SPC procedure. We repeat this simulation 100
times. The boxplots of the achieved false alarm rates from these 100
simulations for different SPC procedures are shown in Figure~\ref{fig:QC0}(a).

As we can see from the plot, the parametric procedure can achieve the
nominal false alarm rate as expected, since the normality assumption is
satisfied in this case. In contrast, the achieved false alarm rate for
the $D_n$ based procedure is far higher than the target value 0.0027.
It is not surprising since all the $\Y_i$ outside the convex hull of
the $\X_i$ will have zero $D_n$ and will be labeled as out-of-control,
but some of those $\Y_i$ may have nonzero $D$ and may have been labeled
as in-control if $D$ was used.
From the plot, we can see that our $R_n$ based procedure can
successfully correct the inflated false alarm rate of the $D_n$ based
procedure and yields the false alarm rate near the target value 0.0027.

We run the same simulations as above on the data generated from the
bivariate elliptical distribution of Section~\ref{sec2.3}. Since the
bivariate elliptical distribution satisfies the conditions of
Theorem~\ref{2}, here we use $R_n$ based on (\ref{eq:md}). Figure~\ref{fig:QC0}(b) shows the corresponding boxplots of the achieved false
alarm rates from 100 simulations for different SPC procedures. As seen
from the plot, the parametric procedure can no longer achieve the
nominal false alarm rate since the normality assumption does not hold
in this case. The $D_n$ based procedure still yields a far higher false
alarm rate than the nominal level, while our $R_n$ based procedure can
achieve the nominal false alarm rate as expected.

%

\begin{figure}

\includegraphics{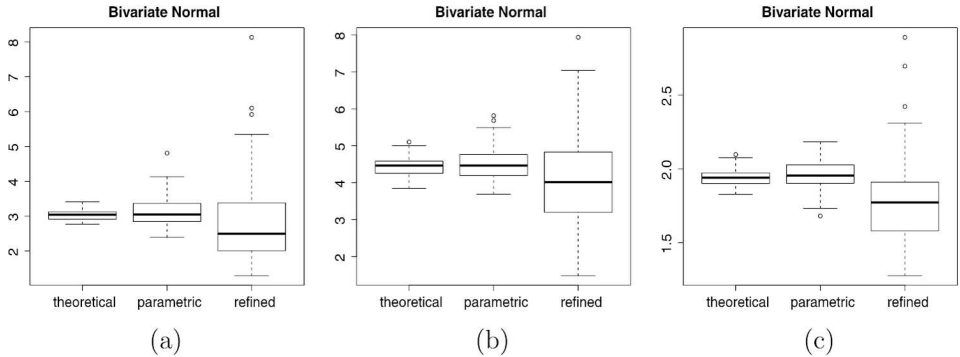}

\caption{The achieved ARLs for the $D$ based procedure, the parametric
procedure and the $R_n$ based procedure under bivariate normal
distribution for \textup{(a)} location change from $(0,0)$ to $(2,2)$; \textup{(b)} scale
increase from 1 to 2; \textup{(c)} both changes in \textup{(a)} and \textup{(b)}.}
\label{fig:QC2}
\end{figure}

To demonstrate the detection power of our $R_n$ based procedure for
process changes, we also carry out the following simulations. Similar
to the above false alarm rate study, we first generate $n=500$
historical data $\X_i$ from the standard bivariate normal distribution.
We then generate $5000$ future observations $\Y_i$ from another
bivariate normal distribution mimicking the following three process
changes: (i) location change from $(0,0)$ to $(2,2)$; (ii) scale increase
from 1 to 2; (iii) both changes in (i) and (ii). Since the $D_n$ based
procedure fails to achieve the nominal false alarm rate, we only
compare the detection power of the parametric procedure and our $R_n$
based procedure. To benchmark the performance, we also include the
procedure based on the theoretical $D$ ($D$ based procedure) in the
comparison. In SPC, a common way to measure the detection power of SPC
procedures is through the average run length (ARL). ARL is the expected
number of times a process needs to be sampled until a specified change
in the process is detected as
out-of-control by the control chart in use.
Figure~\ref{fig:QC2} shows the boxplots of the ARLs from 100
simulations for the three procedures under the three process changes.
As we can see from the plots, the parametric procedure and the $D$
based procedure perform very similarly. Our $R_n$ based procedure
yields slightly smaller ARLs than the $D$ based procedure. This can be
explained by $R_n$'s slightly larger false alarm rate than the nominal
one in Figure~\ref{fig:QC0}(a).

%
%

\begin{figure}

\includegraphics{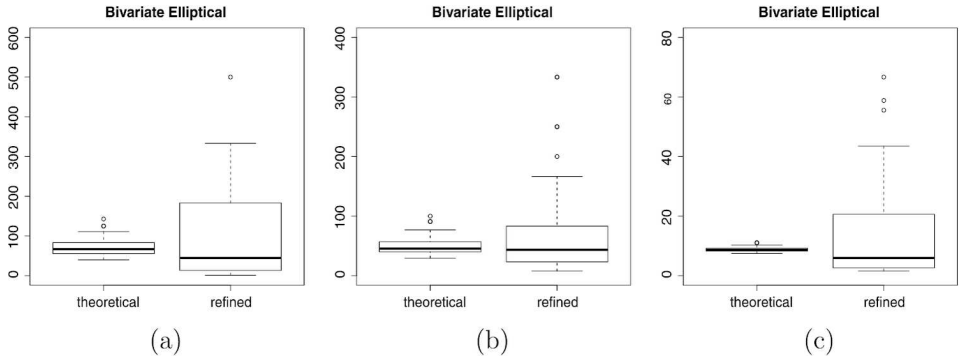}

\caption{The achieved ARLs for the $D$ based procedure and the $R_n$
based procedure under bivariate elliptical distribution for \textup{(a)}
location change from $(0,0)$ to $(4,4)$; \textup{(b)} scale increase from 1 to 2; \textup{(c)}~both changes in \textup{(a)} and \textup{(b)}.}
\label{fig:QC1}
\end{figure}

We repeat the above ARL study on the data generated from the bivariate
elliptical distribution. Similarly, we consider the following three
process changes: (i) location change from $(0,0)$ to $(4,4)$; (ii) scale
increase from 1 to 2; (iii) both changes in (i) and (ii). Since the
parametric procedure does not achieve the nominal false alarm rate in
this bivariate elliptical setting, we only compare the ARLs of the $D$
based procedure and our $R_n$ based procedure. Figure~\ref{fig:QC1}
shows the boxplots of ARLs of the two procedures under different
process changes. As expected, our $R_n$ based procedure performs well
compared with the impractical procedure based on the unknown $D$.

%

\subsection{Classification}\label{sec3.2}
Another application in which the refined half-space depth $R_n$ helps
improve the performance is the classification problem. Classification
is one of the most practical subjects in statistics. It has many
important applications in different fields. For simplicity, we only
focus on two-class classification problem here. In this case, we
observe two training samples $\{\X_1,\ldots,\X_m\}$
and $\{\Y_1,\ldots,\Y_n\}$ from distributions $F$ and $G$,
respectively. The goal of the classification problem is to assign the
future observation $Z$ to either $F$ or $G$ based on some
classification rule built on the two training samples. Recently \citet{LiCueLiu12} developed a nonparametric classification
procedure, called \textit{DD}-classifier, using the \textit{DD}-plot
(depth vs. depth plot) introduced in \citet{LiuParSin99}.
For any two samples, the \textit{DD}-plot plots the depth values of
those pooled sample points with respect to one sample against their
depth values with respect to the other sample. The basic idea behind
the \textit{DD}-classifier is to look for a curve that best separates
the two samples in their \textit{DD}-plot. Since the best separating
curve in the \textit{DD}-classifier is required to pass through the
origin in the \textit{DD}-plot, any future observations having zero
depth values with respect to both samples will be on the separating
curve, indicating that they can be from either sample. Therefore, those
observations will be randomly assigned to either sample. When the $D_n$
of the half-space depth is used in constructing the \textit{DD}-plot,
any point which lies outside of the convex of both samples will have
zero half-space depths with respect to both samples. Based on the
\textit{DD}-classifier, those points will be randomly assigned to
either of the two samples, which will yield roughly a $50\%$
misclassification rate for those points. This simply implies that when
using $D_n$ in the \textit{DD}-classifier one loses all the information
contained in those points. Next, we present a simulation study showing
that the misclassification rate of those points can be improved by
using $R_n$ instead of $D_n$ in the \textit{DD}-classifier.

\begin{figure}

\includegraphics{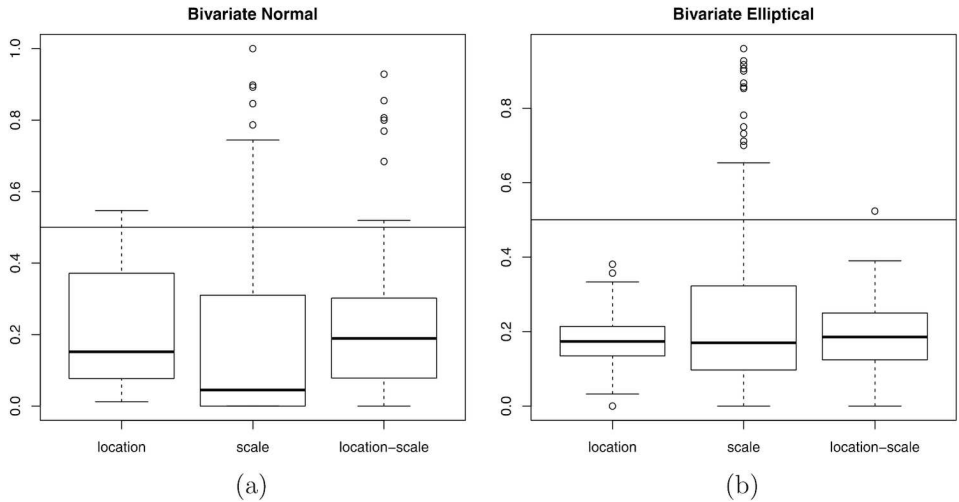}

\caption{The misclassification rate based on $R_n$ under \textup{(a)} bivariate
normal distribution; \textup{(b)} bivariate elliptical distribution.}
\label{fig:cl}
\end{figure}

The first simulation setting we consider is when both $F$ and $G$ are
bivariate normal distributions. We set $F$ as the standard bivariate
normal distribution, and $G$ is another bivariate normal distribution
which differs from $F$ in (i) location; (ii)~scale; (iii) both location
and scale.
(The location difference is 2 for both coordinates; the scale
difference is also 2 for both coordinates.)
For each of the three choices of $G$, we generate a training set
consisting of $m=500$ and $n=500$ observations from $F$ and $G$,
respectively. Based on this training set, we obtain the linear \textit
{DD}-classifier using $R_n$ to construct the \textit{DD}-plot. Another
5000 test observations (2500 from each group) are then generated. Among
those 5000 observations, the misclassification rate for the points
which have zero $D_n$ values with respect to both training samples are
computed. This experiment is repeated 100 times and the
misclassification rates for those points are then summarized in a
boxplot for each choice of $G$ in Figure~\ref{fig:cl}(a).

%

We repeat this simulation on the data where both $F$ and $G$ are
bivariate elliptical distributions; $F$ corresponds to the elliptical
density of Section~\ref{sec2.3}. Again three kinds of differences are
considered: (i) $F$ and $G$ differ in location; (ii) $F$ and $G$ differ
in scale; (iii) $F$ and $G$ differ in both location and scale. (The
location difference is 4 for both coordinates; the scale difference is
2 for both coordinates.) The boxplots of the misclassification rates
for the test observations which have zero $D_n$ values with respect to
both training samples are shown in Figure~\ref{fig:cl}(b).

As mentioned earlier, if $D_n$ is used in the \textit{DD}-classifier,
the misclassification rate for the points which lie outside of the
convex hull of both samples is roughly $50\%$. Therefore, as seen from
Figure~\ref{fig:cl}, the \textit{DD}-classifier paired with $R_n$
substantially improves the classification results for those points.

\section{Concluding remarks}\label{sec4}

We have seen that both applications of the half-space depth in SPC and
classification gain substantially from the proposed refinement $R_n$.
In general, we can expect similar gains from using $R_n$ in statistical
inference methods involving depth ranks or extreme depth contours,
for example, determining $p$-values using depth in \citet{LiuSin97};
constructing multivariate spacings and tolerance regions in \citet{LiLiu08}.

There are many other well-known depth functions [e.g., the spatial
depth \citet{Cha96}, the Mahalanobis depth
\citet{Mah36}, the projection depth \citet{Zuo03}, etc.] which are
not computed from the empirical distribution function, and hence they
do not have the said problem in this paper. While these depths are
useful for many applications, they are either parametric in nature or
lack of the needed distributional properties to ensure the desired
probability masses associated with the central regions formed by the
depth ranks or contours. When these properties are essential, the
applications may be better served by using the two geometric depths.
Case in point are the examples mentioned in the preceding paragraph.
This in part explains the importance in refining the empirical
half-space depth.

It is easy to see that the problem we faced in this paper stems from
the use of the empirical distribution in computing the half-space
probabilities. A natural solution then would be to consider instead a
smoothed version of the empirical distribution that does not have point
masses and is supported on the entire $\mathbb{R}^d$. It is worth
noting that our proposed refinement is in fact such a smoothed version
of the empirical distribution function in the tail, with the smoothing
done by way of extreme value statistics. This extreme-value-theory
based smoothing not only has the advantages of both breaking ties in
the tail and yielding positive values, but, most importantly, it also
produces a statistically much better estimator of the half-space depth
in the tail, as shown in our theorems and applications.

It would be worthwhile to investigate whether the extreme-value-theory
approach proposed in this paper can be modified to refine the empirical
simplicial depth or other depth functions that also use the empirical
counts based on the data. The modifications, if any, would seem quite
nontrivial, since those depth functions do not have such a clear form
of univariate projections as that of the half-space depth.


\section{Proofs}\label{sec5}
\mbox{}
\begin{pf*}{Proof of Theorem~\ref{th:1d}} Write $F^{-1}$ for the
quantile function, the left-continuous inverse of $F$.
We split the region over which the supremum is taken into three
regions: $[F^{-1}(\delta_n), X_{k+1:n}]$, $(X_{k+1:n}, X_{n-k:n})$,
and $[X_{n-k:n},\break F^{-1}((1-\delta_n)^+)]$. Because of symmetry, the
first and last region can be dealt with similarly. Therefore, we
only consider the latter two regions.

For $x \in(X_{k+1:n}, X_{n-k:n})$, we easily see that
%
\begin{eqnarray}
\label{mm} \min \biggl(\frac{F_n(x)}{F(x)},\frac{ S_n(x)}{ S(x)} \biggr)&\leq&
\frac
{\min(F_n(x), S_n(x))}{\min(F(x), S(x))}=\frac{R_n(x)}{D(x)}
\nonumber
\\[-8pt]
\\[-8pt]
\nonumber
&\leq&\max \biggl(\frac{F_n(x)}{F(x)},
\frac{ S_n(x)}{ S(x)} \biggr).
\end{eqnarray}
We have that
%
\begin{equation}
\label{sw}\sup_{x  :  F(x)\geq{k}/{(2n)}}\biggl\llvert \frac{F_n(x)}{F(x)}-1
\biggr\rrvert \stackrel{p} {\to} 0\quad \mbox{and}\quad \sup_{ x  :   S(x)\geq
{k}/{(2n)}}\biggl
\llvert \frac{S_n(x)}{S(x)}-1\biggr\rrvert \stackrel{p} {\to} 0;
\end{equation}
see, for example, Shorack and Wellner
[(\citeyear{ShoWel86}), page~424]. Since $F(X_{k+1:n})> \frac{k}{2n}$ and
$F(X_{n-k:n})< 1- \frac{k}{2n}$ with probability tending to one
($n \to\infty$), it follows from (\ref{sw}) and (\ref{mm})
that
\[
\sup_{X_{k+1:n}<x< X_{n-k:n}}\biggl\llvert \frac{R_n(x)}{D(x)}-1\biggr\rrvert
\stackrel {p} {\to} 0.
\]

Hence, it remains to consider the supremum over the region
$[X_{n-k:n}, F^{-1}((1-\delta_n)+)]$. We have with probability
tending to one, as $n \to\infty$,
\begin{eqnarray*}
&&\sup_{X_{n-k:n}\leq x \leq F^{-1}((1-\delta_n)+)}\biggl\llvert \frac
{R_n(x)}{D(x)}-1\biggr\rrvert
\\
&&\qquad\leq\sup_{\delta_n\leq S(x)\leq2k/n,  S(x)\neq k/n}\biggl\llvert \frac{k}{nS(x)} \biggl(1+\hat
\gamma\frac{x-\hat b}{\hat a} \biggr)^{-1/\hat\gamma
}-1\biggr\rrvert\\
&&\qquad\quad{} +\biggl\llvert
\biggl(1+\hat\gamma\frac{b-\hat b}{\hat a} \biggr)^{-1/\hat\gamma}-1\biggr\rrvert .
\end{eqnarray*}
Write $B_n=\sqrt{k}(\hat b- b)/a$. Then we have $B_n=O_p(1)$; see, for example,
Theorem~2.4.1 in \citet{deHFer06}.
Therefore, to complete the proof of this theorem it suffices to show
%
\begin{equation}
\label{compl}\sup_{\delta_n\leq S(x)\leq2k/n,
S(x)\neq k/n}\biggl\llvert \frac{k}{nS(x)}
\biggl(1+\hat\gamma\frac{x-\hat
b}{\hat a} \biggr)^{-1/\hat\gamma}-1\biggr\rrvert
\stackrel{p} {\to} 0.
\end{equation}

First, we consider the case $\gamma\neq0$. Write $Y_n=\frac{\hat\gamma
}{\gamma}\frac{a}{\hat a}  $. Also, set
\begin{eqnarray}
s=\frac{({(x-b)}/{a})({\gamma}/{(d_n^\gamma-1)})-1}{A}
\nonumber
\\
\eqntext{\mbox {with }\displaystyle d_n=d_n(x)=
\frac{k}{nS(x)} \mbox{ and } A=A(n/k).}
\end{eqnarray}
We have
\begin{eqnarray*}
&&\frac{k}{nS(x)} \biggl(1+\hat\gamma\frac{x-\hat b}{\hat a} \biggr)^{-1/\hat\gamma}
\\
&&\qquad=d_n \biggl(1+Y_n \biggl[\frac{x- b}{a}\gamma-
\frac{\hat b- b}{ a}\gamma \biggr] \biggr)^{-1/\hat\gamma}
\\
&&\qquad=d_n \biggl(1+Y_n \biggl[(1+sA)
\bigl(d_n^\gamma-1\bigr)-\frac{\hat b- b}{ a}\gamma \biggr]
\biggr)^{-1/\hat\gamma}
\\
&&\qquad= \biggl(d_n^{-\hat\gamma}+Y_nd_n^{-\hat\gamma}(1+sA)
\bigl(d_n^\gamma -1\bigr)-\frac{\hat b- b}{ a}\gamma
Y_n d_n^{-\hat\gamma} \biggr)^{-1/\hat
\gamma}
\\
&&\qquad= \biggl[d_n^{\gamma-\hat\gamma} \biggl(d_n^{-\gamma}
\biggl[1-Y_n(1+sA)-\frac{\hat b- b}{ a}\gamma Y_n
\biggr]+Y_n(1+sA) \biggr) \biggr]^{-1/\hat
\gamma}
\\
&&\qquad=: \bigl[ T_1(T_2+T_3)
\bigr]^{-1/\hat\gamma}.
\end{eqnarray*}
We will now prove that $T_1\stackrel{p}{\to}1,
T_2\stackrel{p}{\to}0,    T_3\stackrel{p}{\to}1$, all
\textit{uniformly for $x$ such that $\delta_n\leq S(x)\leq
2k/n$} [$S(x)\neq k/n$]. This will yield (\ref{compl}) for
$\gamma\neq0$.

We have
\[
T_1=d_n^{\gamma-\hat\gamma}=d_n^{-\Gamma_n/\sqrt{k}}=
\exp \biggl( \frac
{-\Gamma_n}{\sqrt{k}}\log\frac{k}{nS(x)} \biggr).
\]
Observe that
\[
\biggl\llvert \frac{-\Gamma_n}{\sqrt{k}}\log\frac{k}{nS(x)}\biggr\rrvert \leq
\frac
{|\Gamma_n|}{\sqrt{k}}\biggl\llvert \log\frac{k}{nS(x)}\biggr\rrvert \stackrel{p} {
\to}0.
\]
Hence, $T_1\stackrel{p}{\to}1$.
Consider $T_3=Y_n(1+sA)$.
We have $Y_n\stackrel{p}{\to}1$ and $A(n/k)\to0$. Hence, (\ref
{finite}) yields $T_3\stackrel{p}{\to}1$.
Finally,
\begin{eqnarray*}
T_2&=&\frac{d_n^{-\gamma}}{\sqrt{k}}\sqrt{k} \biggl(1-Y_n(1+sA)-
\frac
{B_n}{\sqrt{k}} \gamma Y_n \biggr)
\\
&=& \frac{d_n^{-\gamma}}{\sqrt{k}} \biggl(\sqrt{k} \biggl[1- \biggl(1+O_p
\biggl(\frac{1}{\sqrt{k}} \biggr) \biggr) \biggl(1+O \biggl(\frac{1}{\sqrt{k}}
\biggr) \biggr) \biggr]-B_n \gamma Y_n \biggr)
\\
&=& \frac{(nS(x))^\gamma}{k^{\gamma+1/2}}O_p(1)=o_p(1).
\end{eqnarray*}

Consider now the case $\gamma=0$. By convention $(d_n^\gamma-1)/\gamma
=\log d_n$ now. Write
\begin{eqnarray*}
Q&:=&\frac{k}{nS(x)} \biggl(1+\hat\gamma\frac{x-\hat b}{\hat a}
\biggr)^{-1/\hat\gamma}
\\
&=&d_n \biggl(1+\frac{a}{\hat a} \biggl[\frac{x- b}{a}\hat
\gamma-\frac
{\hat b- b}{ a}\hat\gamma \biggr] \biggr)^{-1/\hat\gamma}
\\
&=&d_n \biggl(1+\hat\gamma\frac{a}{\hat a}(\log d_n)
(1+sA)-\hat\gamma \frac{B_n}{ \sqrt{k}} \frac{a}{\hat a} \biggr)^{-1/\hat\gamma}.
\end{eqnarray*}
Hence,
\[
\log Q=\log d_n-\frac{1}{\hat\gamma} \log \biggl(1+\hat\gamma
\frac{a}{\hat a}(\log d_n) (1+sA)-\hat \gamma\frac{B_n}{ \sqrt{k}}
\frac{a}{\hat a} \biggr).
\]
We obtain
\[
|\log Q|=\log \biggl( \frac{k}{n\delta_n} \biggr) \cdot O_p \biggl(
\frac
{1}{\sqrt{k}} \biggr)+\log^2 \biggl( \frac{k}{n\delta_n} \biggr)
\cdot O_p \biggl(\frac{1}{\sqrt{k}} \biggr)\stackrel{p} {\to}0.
\]
Hence,
$Q\stackrel{p}{\to}1$, uniformly for $x$ such that
$\delta_n\leq S(x)\leq2k/n$ ($S(x)\neq k/n$). This proves
(\ref{compl}) for $\gamma=0$.
\end{pf*}

For the proof of Theorem~\ref{2}, we need two lemmas. In the sequel, we
assume that the conditions of Theorem~\ref{2} are in force. Write $\Theta=\{\u
\in\mathbb{R}^d: \Vert \u\Vert =1\}$ for the unit sphere.

\begin{lemma}\label{l1} For all $r>0$,
\[
\lim_{t\to\infty} \sup_{\u\in\Theta} \biggl\llvert
\frac{\mathbb{P}(\X
_1\in tH_{r,\u})}{t^{-\alpha}}-c \nu(H_{r,\u})\biggr\rrvert =0.
\]
\end{lemma}

\begin{pf}
Fix $r>0$. Combining (\ref{nu}) and
(\ref{radi}) we have that for all $\u\in\Theta$,
%
\begin{equation}
\label{34}\lim_{t\to\infty} \frac{\mathbb{P}(\X_1\in
tH_{r,\u})}{t^{-\alpha}}=c
\nu(H_{r,\u}).
\end{equation}
Assume this convergence does not hold uniformly
in $\u\in\Theta$. Then there exist sequences $\u_m\to\vv$ and $t_m\to
\infty$ such that
%
\begin{equation}
\label{nocon}\mathbb{P}(\X_1\in t_mH_{r,\u
_m})/t_m^{-\alpha}
\mbox{ does \textit{not} converge to } c \nu(H_{r,\vv
}), \mbox{ as } m \to
\infty.
\end{equation}
W.l.o.g. we assume that $\vv=(1, 0, \ldots, 0)$.

We show that (\ref{nocon}) cannot hold by showing, for $\u\in\Theta$,
%
\begin{equation}
\label{ep} 
\frac{\mathbb{P}(\u^T\X_1\geq {tr})}{\mathbb{P}(X_{1,1}\geq {tr})} \to1\qquad\mbox{if }
u_1\to1, t\to\infty.
\end{equation}
Because if the latter convergence holds, then if $u_1\to1, m\to\infty$,
%
\begin{equation}
\label{useful} \frac{\mathbb{P}(\u^T\X_1\geq t_mr)}{t_m^{-\alpha}} =\frac{\mathbb{P}(\u^T\X_1\geq t_mr)}{
\mathbb{P}(X_{1,1}\geq t_mr)} \cdot \frac{\mathbb{P}(X_{1,1}\geq
t_mr)}{t_m^{-\alpha}}
\to1\cdot c\nu(H_{r,\vv}).
\end{equation}

Hence, it remains to show (\ref{ep}). Write $\varepsilon=1-u_1$. Then
$\varepsilon\to0$.
We have
\begin{eqnarray*}
&&\mathbb{P}\bigl(\u^T\X_1\geq {tr}\bigr)
\\
&&\qquad= \mathbb{P}\bigl(\u^T\X_1\geq {tr}, X_{1,1} <
\bigl(1-\varepsilon ^{1/4}\bigr){tr}\bigr)+\mathbb{P}\bigl(
\u^T\X_1\geq {tr}, X_{1,1} \geq\bigl(1-\varepsilon
^{1/4}\bigr){tr}\bigr)
\\
&&\qquad\leq \mathbb{P}\bigl(\u^T\X_1-(1-
\varepsilon)X_{1,1} \geq\varepsilon^{1/4}{tr}\bigr) +\mathbb{P}
\bigl( X_{1,1} \geq\bigl(1-\varepsilon^{1/4}\bigr){tr}\bigr)
\\
&&\qquad\leq \sum_{j=2}^d \mathbb{P}
\bigl(u_jX_{1,j}\geq\varepsilon ^{1/4}{tr}/(d-1)
\bigr)+\mathbb{P}\bigl( X_{1,1} \geq\bigl(1-\varepsilon^{1/4}
\bigr){tr}\bigr)
\\
&&\qquad\leq \sum_{j=2}^d \mathbb{P}
\bigl(|X_{1,j}|\geq\varepsilon^{-1/4}{tr}/\bigl(\sqrt {2}(d-1)
\bigr)\bigr)+\mathbb{P}\bigl( X_{1,1} \geq\bigl(1-\varepsilon^{1/4}
\bigr){tr}\bigr).
\end{eqnarray*}
Hence,
\begin{eqnarray*}
&&\frac{\mathbb{P}(\u^T\X_1\geq {tr})}{\mathbb{P}(X_{1,1}\geq {tr})}
\\
&&\qquad\leq \frac{\sum_{j=2}^d \mathbb{P}(|X_{1,j}|\geq\varepsilon
^{-1/4}{tr}/(\sqrt{2}(d-1)))}{t^{-\alpha}}\cdot\frac{t^{-\alpha}}{
\mathbb{P}(X_{1,1}\geq {tr})}
\\
&&\qquad\quad{} +\frac{\mathbb
{P}(X_{1,1}\geq(1-\varepsilon^{1/4}){tr})}{t^{-\alpha}}\cdot\frac
{t^{-\alpha}}{\mathbb{P}(X_{1,1}\geq {tr})}
\\
&&\qquad\to 0\cdot\frac{1}{c \nu(H_{r,\vv})}+ c \nu(H_{r,\vv})\cdot
\frac
{1}{c \nu(H_{r,\vv})}=1.
\end{eqnarray*}
Similarly, we have
\begin{eqnarray*}
&&\mathbb{P}\bigl(\u^T\X_1\geq {tr}\bigr)
\\
&&\qquad\geq \mathbb{P}\bigl(X_{1,1}\geq\bigl(1+\varepsilon^{1/4}
\bigr){tr}\bigr)-\mathbb{P}\bigl(\u^T\X _1<
{tr},X_{1,1}\geq\bigl(1+\varepsilon^{1/4}\bigr){tr}\bigr)
\\
&&\qquad\geq \mathbb{P}\bigl(X_{1,1}\geq\bigl(1+\varepsilon^{1/4}
\bigr){tr}\bigr) -\mathbb{P}\bigl(\u ^T\X_1-(1-
\varepsilon)X_{1,1} \leq-\varepsilon^{1/4}{tr}/2\bigr)
\end{eqnarray*}
and
\begin{eqnarray*}
&&\frac{\mathbb{P}(\u^T\X_1\geq {tr})}{\mathbb{P}(X_{1,1}\geq {tr})}
\\
&&\qquad\geq \frac{\mathbb{P}(X_{1,1}\geq(1+\varepsilon^{1/4}){tr})}{\mathbb
{P}(X_{1,1}\geq {tr})} -\frac{\mathbb{P}(\u^T\X_1-(1-\varepsilon)X_{1,1}
\leq-\varepsilon^{1/4}{tr}/2)}{\mathbb{P}(X_{1,1}\geq {tr})}
\\
&&\qquad\to c \nu(H_{r,\vv})\cdot\frac{1}{c \nu(H_{r,\vv})}-0\cdot\frac{1}{c
\nu(H_{r,\vv})}=1.
\end{eqnarray*}
This completes the proof of (\ref{ep}).
\end{pf}

%
%
%

Define the function $g$ by $g(\u)=c\nu(H_{1,\u})$ and let
$V_\u(t)=F_\u^{-1}(1-1/t),  t>1$, be the tail quantile function
corresponding to $F_\u$.
%
\begin{lemma}\label{l3} We have
%
\begin{equation}
\label{U}\lim_{t\to\infty} \sup_{\u\in\Theta} \biggl
\llvert \frac{V_\u(t)}{t^{1/\alpha}}- \bigl(g(\u)\bigr)^{1/\alpha}\biggr\rrvert =0.
\end{equation}
\end{lemma}

\begin{pf}
Lemma~\ref{l1}, with $r=1$, yields
%
\begin{equation}
\label{unif} \lim_{s\to\infty} \sup_{\u\in\Theta} \biggl
\llvert \frac{1-F_\u
(s)}{s^{-\alpha}}- g(\u)\biggr\rrvert =0.
\end{equation}
Observe that
$V_\u(t)=s$ implies $F_\u(s)=1-1/t$.
Hence
%
\begin{equation}
\label{UF}\frac{V_\u(t)}{t^{1/\alpha}}=s\bigl(1-F_{\u
}(s)
\bigr)^{1/\alpha}= \biggl(\frac{1-F_\u(s)}{s^{-\alpha}} \biggr)^{1/\alpha}.
\end{equation}
Also observe that assumption (\ref{inf}) and $\nu(H_{1,\u})\leq1, \u
\in\Theta$, yield
\[
0<\inf_{\u\in\Theta} g(\u)\leq\sup_{\u\in\Theta} g(\u)\leq c
<\infty.
\]
Combining this with (\ref{UF}) and (\ref{unif}) easily yields (\ref
{U}).
\end{pf}

\begin{pf*}{Proof of Theorem~\ref{2}} We will prove that, as $n \to
\infty$,
%
\begin{eqnarray}
\label{i} \sup_{D(\x)\geq{k}/{(2n)}}\biggl\llvert \frac{D_n(\x)}{D(\x)}-1
\biggr\rrvert &\stackrel{p} {\to}&0 \quad\mbox{and}
\\
\label{ii} \sup_{R_n(\x)<{k}/{n},  D(\x)\geq\delta_n}\biggl\llvert \frac{R_n(\x
)}{D(\x)}-1
\biggr\rrvert& \stackrel{p} {\to}&0.
\end{eqnarray}
To show that (\ref{i}) and (\ref{ii}) imply (\ref{theo2}), it is
sufficient to show that (\ref{i}) implies
%
\begin{equation}
\label{c} \sup_{D_n(\x)\geq{k}/{n}}\biggl\llvert \frac{D_n(\x)}{D(\x)}-1
\biggr\rrvert \stackrel{p} {\to}0
\end{equation}
and to recall that if $D_n(\x)\geq k/n$ or $R_n(\x)\geq k/n$, then
$D_n(\x)=R_n(\x)$.

Assume (\ref{i}) holds.
It follows from \citet{DonGas92}, that\break  $\sup_\x D_n(\x)\geq
1/(d+1)$, with probability 1.
Hence, for large $n$, any point $\hat\x$ with maximum depth $D_n$,
satisfies $D_n(\hat\x)\geq k/n$ and, with probability tending to one,
$D(\hat\x)\geq k/n$, because of the uniform consistency of $D_n$. Now
assume for some $\x$, $D_n(\x)\geq k/n$ and $D(\x)<k/(2n)$. Then, with
probability tending to one, we can find $\x_0$ on the straight line
connecting $\hat\x$ and $\x$, such that $D(\x_0)=k/(2n)$ and because
of (\ref{i}), $D_n(\x_0)\leq3k/(4n)$. It is well known that $D_n$ has
the ``monotonicity relative to deepest point'' property
[see, e.g., \citet{ZuoSer00}], and hence $D_n(\x)\leq D_n(\x_0)\leq3k/(4n)$.
Contradiction. Hence (\ref{c}).

It remains to prove (\ref{i}) and (\ref{ii}). We begin with (\ref{i}).
First, we show that
%
\begin{equation}
\label{alex} P \Bigl(\bigcup \bigl\{H:P(H)\leq s\bigr\} \Bigr)=O(s)\qquad \mbox{as } s
\downarrow0.
\end{equation}
Define $r_0=(c\inf_{\u\in\Theta}\nu(H_{1,\u})/2)^{1/\alpha}$. Lemma~\ref{l1}
yields that, uniformly in $\u\in\Theta$,
\[
\lim_{s\downarrow0} \frac{\mathbb{P}(\X_1\in s^{-1/\alpha}H_{r_0,\u
})}{s}=c \nu(H_{r_0,\u})=cr_0^{-\alpha}
\nu(H_{1,\u})\geq2.
\]
Hence, for small enough $s$ and uniformly in $\u\in\Theta$,
$\mathbb{P}(\X_1\in s^{-1/\alpha}H_{r_0,\u})>s$. For $\u\in\Theta$,
let $r_1$ be the smallest $r$ such that
$P(H_{r,\u})=s$. Then for small enough $s$, $H_{r_1,\u}\subset
s^{-1/\alpha}H_{r_0,\u}$.
Hence, by (\ref{nu}) and (\ref{radi}),
\begin{eqnarray*}
\frac{P (\bigcup_{\u\in\Theta}H_{r_1,\u} )}{s}&\leq&\frac{P
(\bigcup_{\u\in\Theta}s^{-1/\alpha}H_{r_0,\u} )}{s} =\frac{P (s^{-1/\alpha}
\bigcup_{\u\in\Theta}H_{r_0,\u}
)}{s}\\
&\to& c\nu\biggl(
\bigcup_{\u\in\Theta}H_{r_0,\u}\biggr)<\infty,\qquad s
\downarrow0,
\end{eqnarray*}
which implies (\ref{alex}).

Using (\ref{alex}), we obtain from Theorem~5.1 in \citet{Ale87},
with the $\gamma_n$ there equal to $k/(2n)$, that
%
\begin{equation}
\label{alexth} \sup_{H:  P(H)\geq{k}/{(2n)}}\biggl\llvert \frac{P_n(H)}{P(H)}-1
\biggr\rrvert \stackrel{p} {\to} 0\qquad \mbox{as } n \to\infty.
\end{equation}

Denote with $H_\x$ a half-space with $\x$ on its boundary. We have
\[
\frac{D_n(\x)}{D(\x)}=\frac{\inf_{H_\x} P_n(H_\x)}{\inf_{H_\x} P(H_\x
)}\geq\inf_{H_\x}
\frac{ P_n(H_\x)}{ P(H_\x)}
\]
and, with $\varepsilon>0$, for some $H_\x $,
\[
\frac{D_n(\x)}{D(\x)}\leq(1+\varepsilon) \frac{ P_n(H_\x)}{ P(H_\x)}.
\]
This, in combination with (\ref{alexth}), yields (\ref{i}).

Finally, we consider (\ref{ii}). Write $p_\u(w)=\mathbb{P}(\u^T\X_1\geq
w)=P(H_{w,\u})$. We first show
%
\begin{equation}
\label{pdp} \sup_{w, \u\in\Theta:  \delta_n\leq p_\u(w)\leq{2k}/{n}}\biggl\llvert \frac{p_{n,\u}(w)}{p_\u(w)}-1
\biggr\rrvert \stackrel{p} {\to} 0\qquad\mbox {as } n \to\infty.
\end{equation}
Write $d_\u(w)=k/(np_\u(w))$.
Then
%
\begin{eqnarray}
\label{pdp2} \frac{p_{n,\u}(w)}{p_\u(w)}&=&\frac{k}{np_\u(w)} \biggl(\frac{w}{V_\u(
{n}/{k})}
\frac{V_\u({n}/{k})}{W_{n-k:n}} \biggr)^{-\widehat\alpha}
\nonumber
\\[-8pt]
\\[-8pt]
\nonumber
& = &\biggl(d_\u^{-1/\widehat\alpha}(w)
\frac{w}{V_\u({n}/{k})} \biggr)^{-\widehat\alpha} \biggl(\frac{V_\u({n}/{k})}{W_{n-k:n}}
\biggr)^{-\widehat\alpha}.
\end{eqnarray}
%
It follows from Lemmas \ref{l1} and \ref{l3} that
\[
\lim_{n\to\infty} \sup_{w, \u\in\Theta:  \delta_n\leq p_\u(w)\leq
{2k}/{n}} \biggl\llvert
d_\u^{-1/\alpha}(w)\frac{w}{V_\u({n}/{k})}-1\biggr\rrvert =0.
\]
Using this, (\ref{gamma}) and $\log(n\delta_n)/\sqrt{k}\to0$, we obtain
%
\begin{equation}
\label{ratio} \sup_{w,  \u\in\Theta:  \delta_n\leq p_\u(w)\leq{2k}/{n}}\biggl\llvert \biggl(d_\u^{-1/\widehat\alpha}(w)
\frac{w}{V_\u({n}/{k})} \biggr)^{-\widehat\alpha} -1\biggr\rrvert \stackrel{p} {\to} 0,\qquad
\mbox{as } n \to\infty.
\end{equation}
Denote with $G_{\u,n}$ the empirical distribution function of the
uniform-$(0,1)$ random variables $F_u(\u^T\X_i), i=1,\ldots, n$, and with
$G^{-1}_{\u,n}$ the corresponding quantile function. It follows from
(\ref{alexth}) by routine arguments that
\[
\sup_{\u\in\Theta} \biggl\llvert \frac{1-G^{-1}_{\u,n}(1-k/n)}{k/n}-1\biggr\rrvert
\stackrel{p} {\to} 0,
\]
and hence, by Lemma~\ref{l3} and (\ref{gamma}), that
%
\begin{equation}
\label{g} \sup_{\u\in\Theta} \biggl\llvert \biggl(
\frac{V_\u(
{n}/{k})}{W_{n-k:n}} \biggr)^{-\widehat\alpha}-1\biggr\rrvert \stackrel{p} {\to} 0.
\end{equation}
Combination of (\ref{pdp2}), (\ref{ratio}) and (\ref{g}), yields (\ref{pdp}).

Now we turn to (\ref{ii}). Let $\x$ be such that $R_n(\x)<k/n$. Then
%
\begin{eqnarray}
\label{geq} \frac{R_n(\x)}{D(\x)}&=&\frac{\inf_{\u\in\Theta:  1-\widehat F_{\u}
(\u^T\x^-)<k/n}1-\widehat F_{\u} (\u^T\x^-)}{\inf_{\u\in\Theta} p_{\u
}(\u^T\x)}
\nonumber
\\[-8pt]
\\[-8pt]
\nonumber
&\geq&\frac{\inf_{\u\in\Theta:  p_{n,\u}(\u^T\x)<k/n}   p_{n,\u}(\u
^T\x)}{\inf_{\u\in\Theta:  p_{n,\u}(\u^T\x)<k/n}    p_{\u}(\u^T\x)} \geq\inf_{\u\in\Theta:  p_{n,\u}(\u^T\x)<k/n}
\frac{p_{n,\u}(\u^T\x
)}{p_{\u}(\u^T\x)}.
\end{eqnarray}
Next, we show
that with probability tending to one ($n\to\infty$),
%
\begin{equation}
\label{z}\inf_{\u\in\Theta:  p_{n,\u}(\u^T\x
)<k/n}\frac{p_{n,\u}(\u^T\x)}{p_{\u}(\u^T\x)}\geq\inf
_{\u\in\Theta
:  p_{\u}(\u^T\x)\leq2k/n}\frac{p_{n,\u}(\u^T\x)}{p_{\u}(\u^T\x)}.
\end{equation}
Assume for some $\x$ and $\u\in\Theta$, $p_{n,\u}(\u^T\x)<k/n$ and
$p_{\u}(\u^T\x)> 2k/n$. Then there exists an $\x_0$ of the form $\x
+\tilde c\u$, for some $\tilde c>0$, such that $p_{\u}(\u^T\x_0)=
2k/n$. Hence, with probability tending to one because of (\ref{pdp}),
$p_{n,\u}(\u^T\x_0)\geq3k/(2n)$ and, therefore,
$p_{n,\u}(\u^T\x)\geq3k/(2n)$. Contradiction. Hence, we have (\ref
{z}). Combining (\ref{z}) with (\ref{geq}) and (\ref{pdp}), yields
%
\begin{equation}
\label{max}\sup_{R_n(\x)<{k}/{n},  D(\x)\geq\delta
_n} \biggl(1-\frac{R_n(\x)}{D(\x)} \biggr)
\vee 0\stackrel{p} {\to}0.
\end{equation}

Let $\varepsilon\in(0,1)$ and let $\x$ be such that $R_n(\x)<k/n$ and
$D(\x)\geq\delta_n$. We have for some $\u_0$
that
\[
\frac{R_n(\x)}{D(\x)}\leq(1+\varepsilon/2) \frac{R_n(x)}{p_{\u_0}(\u
_0^T\x)} \leq(1+
\varepsilon/2) \frac{1-\widehat F_{\u_0}
(\u_0^T\x^-)}{p_{\u_0}(\u_0^T\x)} ,
\]
with
$p_{\u_0}(\u_0^T\x)\geq\delta_n$. If $p_{\u_0}(\u_0^T\x)\leq
k/(2n)$, then with probability tending to one, (\ref{pdp})
yields that $1-\widehat F_{\u_0}
(\u_0^T\x^-)=p_{n,\u_0}(\u_0^T\x)$, and hence that
$R_n(\x)/\break D(\x)\leq1+\varepsilon$. In case
$p_{\u_0}(\u_0^T\x)>k/(2n)$, we have, using (\ref{max}), that
with probability tending to one that
$k/(2n)<p_{\u_0}(\u_0^T\x)\leq3D(\x)/2\leq2R_n(\x)<2k/n$.
Hence, combining $1-\widehat F_{\u_0} (\u_0^T\x^-)\leq(1- F_{n,
\u_0} (\u_0^T\x^-))\vee p_{n,\u_0}(\u_0^T\x)$ with
(\ref{alexth}) and (\ref{pdp}), we obtain that with
probability tending to one, $R_n(\x)/D(\x)\leq1+\varepsilon$.
Hence, we have shown
%
\[
\sup_{R_n(\x)<{k}/{n},  D(\x)\geq\delta
_n} \biggl(\frac{R_n(\x)}{D(\x)}-1 \biggr) \vee 0
\stackrel{p} {\to}0.
\]
This, in combination with (\ref{max}), yields (\ref{ii}).
\end{pf*}

\section*{Acknowledgements} We thank an Associate Editor and three
referees for many insightful questions and comments, which helped
improve greatly this paper. We also thank Laurens de Haan for his
helpful discussions, in particular about the univariate second-order
condition.


\printaddresses
\end{document}